\documentclass[12pt,a4paper]{article}

\usepackage{amsmath,amsfonts,amssymb}
\usepackage{bbm}
\usepackage{cases}

\usepackage{kotex}
\usepackage{lipsum}
\usepackage{subcaption}

\usepackage{graphicx}
\usepackage{ulem}

\newtheorem{defn}{Definition}[section]
\newtheorem{theo}[defn]{Theorem}

\newtheorem{prop}[defn]{Proposition}

\newtheorem{rem}[defn]{Remark}

\newenvironment{proof}{{\bf Proof }}{{\vskip 0.1cm \hfill$\Box$}}

\begin{document} 

\noindent
{\Large \bf Improved contraction estimates for linear elliptic diffusion-advection-reaction equations}
\\ \\
\bigskip
\noindent
{\bf Haesung Lee}  \\
\noindent
{\bf Abstract.}  
This paper analyzes contraction estimates for the unique strong solution to a linear elliptic diffusion-advection-reaction equation in a bounded open subset of $\mathbb{R}^d$ with $d \geq 3$, where the zero-order coefficient lies in $L^{\frac{d}{2}}$ and the drift coefficient belongs to $L^p$ with $p \in (d, \infty)$. Compared with existing results, we obtain contraction estimates with improved constants in a quantitative sense. In particular, for the $L^\infty$-contraction estimate, the constant is shown to have the explicit form $\frac{1}{\alpha}$, depending only on the constant $\alpha > 0$ in the zero-order term. The key step is to represent the solution as the image of a weighted resolvent associated with a Dirichlet form.


\section{Introduction} 
In this paper, we investigate contraction estimates for the unique strong solution to the following homogeneous Dirichlet boundary value problem for a linear elliptic diffusion-advection-reaction equation in non-divergence form, 
\begin{equation} \label{maineqpa}
\left\{
\begin{alignedat}{2}
- \Delta u + \langle \mathbf{H}, \nabla u \rangle + (c+\alpha) u &= f && \quad \text{in } U, \\ 
u &= 0 && \quad \text{on } \partial U,
\end{alignedat}
\right.
\end{equation}
where $U$ and the coefficients satisfy the following assumption {\bf (A)}. 
\centerline{}
\centerline{}
\noindent
{\bf (A)}
{\it
$U$ is a bounded open subset of $\mathbb{R}^d$, $d \geq 3$, with boundary of class $C^{1,1}$. $\mathbf{H} \in L^p(U, \mathbb{R}^d)$ with $p \in (d, \infty)$ and $c \in L^{\frac{d}{2}}(U)$ satisfies $c \geq 0$ in $U$. $\alpha \in (0, \infty)$ is a constant and $f \in L^q(U) \cap L^{\theta}(U)$ for some $q \in (\frac{d}{2}, \infty)$ and $\theta \in [1, \infty]$.
}
\centerline{}
\centerline{}
Here, $u$ is called a strong solution to \eqref{maineqpa} if 
$u \in H^{2,1}(U) \cap H^{1,2}_0(U)$ and satisfies
\[
- \Delta u(x) + \langle \mathbf{H}(x), \nabla u(x) \rangle + (c(x)+\alpha) u(x)
= f(x) \quad \text{for a.e. } x \in U.
\]
On the other hand, $u$ is called a weak solution to \eqref{maineqpa} if 
$u \in H^{1,2}_0(U)$ satisfies
\[
\int_{U} \langle \nabla u, \nabla \varphi \rangle\,dx
+ \int_{U} \langle \mathbf{H}, \nabla u \rangle \varphi\,dx
+ \int_{U} (c+\alpha)u \varphi\,dx
= \int_{U} f \varphi\,dx
\quad \text{for all } \varphi \in C_0^{\infty}(U).
\]
The main motivation of this paper is to provide a quantitative analysis of the constant $C$
in the following $L^2$-contraction estimate:
\begin{equation*} 
\|u\|_{L^2(U)} \leq \frac{C}{\alpha} \|f\|_{L^2(U)},
\end{equation*}
where $u$ denotes the strong solution to \eqref{maineqpa}.
This is especially relevant in the context of numerical methods based on artificial intelligence (AI), such as physics-informed neural networks (PINNs) (see \cite{MZ22, YL24} and references therein), which have recently attracted considerable attention. In many situations, finding an explicit value of such a constant $C>0$ is highly nontrivial. Indeed, it is frequently obtained only implicitly, relying on uniqueness results together with compactness arguments (cf. \cite[Section 6.2, Theorem 6]{Ev10}, \cite[Lemma 9.17]{GT}).
\\
Before studying the strong solution to \eqref{maineqpa}, it is natural to first explore the corresponding weak solution. In the classical literature (cf. \cite{S65}), for the more general case with diffusion coefficients $A=(a_{ij})_{1 \leq i,j \leq d}$ in place of the identity matrix, the existence and uniqueness of a weak solution to \eqref{maineqpa} are obtained under the assumption that $\alpha$ is sufficiently large. The key method used there is based on the Lax-Milgram theorem, employing the coercivity of the associated bilinear form (see \cite[Chapter 2]{BKRS15}). Later in \cite{T73}, Trudinger considered the case $\alpha=0$ and developed the weak maximum principle (cf. \cite{T77}), from which the existence and uniqueness of weak solutions to \eqref{maineqpa} can be shown. In contrast, a very general result for strong solutions is established in \cite{V92}, making use of properties of the Newtonian potential.\\
More recently, \cite{L24} considers the case where the drift vector field has negative weak divergence with respect to the Lebesgue measure, proving the existence and uniqueness of weak solutions to \eqref{maineqpa}. Subsequently, \cite{L25jm} extends this to more general diffusion coefficients and singular zero-order coefficient $c \in L^1(U)$. In particular, \cite[Theorem 1.2]{L25jm} shows the existence and uniqueness of strong solutions to \eqref{maineqpa}. Precisely, in \cite{L25jm}, a divergence-free transformation (\cite[Theorems 3.1, 3.2]{L25jm}) was developed, enabling the reformulation of \eqref{maineqpa} into a linear elliptic partial differential equation with weak divergence-free drift coefficients (see Theorem \ref{divfretransform}):
\begin{equation} \label{maineq3}
\left\{
\begin{alignedat}{2}
- \operatorname{div}(\rho \nabla u) + \langle \rho \mathbf{B}, \nabla u \rangle + \rho (c+\alpha) u &= \rho f && \quad \text{in } U, \\ 
u &= 0 && \quad \text{on } \partial U,
\end{alignedat}
\right.
\end{equation}
where $\rho \in H^{1,2}(U) \cap C(\overline{U})$ is a strictly positive function and $\mathbf{B}:= \mathbf{H} + \frac{1}{\rho} \nabla \rho$, so that $\rho \mathbf{B}$ is a weak divergence-free vector field. By exploiting \eqref{maineq3}, it was shown in \cite[Theorem 1.1(iii)]{L25jm} that there exists a constant $K_1 \geq 1$ which only depends on $d$, $U$, $p$, $\|\mathbf{H}\|_{L^p(U)}$, such that for each $f \in L^q(U) \cap L^{\theta}(U)$ with $q \in (\frac{d}{2}, \infty)$ and $\theta \in [1, \infty]$,
\begin{equation} \label{eltwoesljm}
\|u \|_{L^{\theta}(U)} \leq \frac{K_1}{\alpha}\|f \|_{L^{\theta}(U)}.
\end{equation}
In \cite{L25jm}, the estimate \eqref{eltwoesljm} is obtained by relying on contraction estimates with respect to the Lebesgue measure $dx$ developed in \cite{L24}, rather than on contraction estimates with respect to a weighted measure $\rho\,dx$.
More precisely, since the zero-order coefficient $\rho (c+\alpha)$ in \eqref{maineq3} is bounded below by $\alpha \min_{\overline{U}} \rho$, we may apply the contraction estimates with respect to the Lebesgue measure $dx$ from \cite{L24} to obtain
\[
\|u\|_{L^{\theta}(U)} \leq \frac{1}{\alpha \min_{\overline{U}} \rho} \|\rho f\|_{L^{\theta}(U)}
\leq \frac{\max_{\overline{U}} \rho}{\alpha \min_{\overline{U}} \rho} \| f\|_{L^{\theta}(U)}.
\]
Moreover, the constant $K_1 \geq 1$ arises from the Harnack inequality
\[
\max_{\overline{U}} \rho \leq K_1 \min_{\overline{U}} \rho 
\]
and therefore the estimate \eqref{eltwoesljm} follows. 
\\
In this paper, we take a different approach by directly identifying the solution to \eqref{maineqpa} as the image of a weighted resolvent, which requires the construction of the associated Dirichlet form (cf. \cite{MR92, FOT11}).  More precisely, under assumption {\bf (A)}, we construct the Dirichlet form $(\mathcal{E}, H^{1,2}_0(U))$ on $L^2(U,\rho\,dx)$ associated with \eqref{maineq3}, defined by
\[
\mathcal{E}(v,w)
= \int_U \langle \nabla v, \nabla w \rangle \rho\,dx
+ \int_U \langle \mathbf{B}, \nabla v \rangle w \rho\,dx
+ \int_U c\, v w  \rho\,dx,
\quad v,w \in H^{1,2}_0(U).
\]
We then construct the associated sub-Markovian resolvent of contractions $(G_\alpha)_{\alpha>0}$ on $L^2(U,\rho\,dx)$, which extends to $L^{\theta}(U,\rho\,dx)$, and identify the solution to \eqref{maineq3} with $G_\alpha f$. This approach, based on contraction estimates with respect to the weighted measure $\rho\,dx$, leads to the main results stated below.
\begin{theo} \label{mainexisthm}
Assume {\bf (A)}.  Then, there exists $u \in H^{2, \frac{d}{2}}(U) \cap H^{1,2}_0(U) \cap L^{\infty}(U)$ such that $u$ is a unique strong solution to \eqref{maineqpa} (cf. \cite[Theorem 1.2]{L25jm}). Moreover, there exists $\mu=\rho dx$ with $\rho \in H^{1,p}(U) \cap C(\overline{U})$ and $\rho(x)>0$ for all $x \in \overline{U}$ such that
\begin{equation} \label{weightedcontra}
\|u\|_{L^{\theta}(U, \mu)} \leq \frac{1}{\alpha} \|f\|_{L^{\theta}(U, \mu)}.
\end{equation}
In particular,
\begin{equation} \label{contrawithrele}
\|u\|_{L^{\theta}(U)} \leq \frac{K_1^{\frac{1}{\theta}}}{\alpha} \|f\|_{L^{\theta}(U)},
\end{equation}
where $K_1 \geq 1$ is a constant as in Theorem \ref{existrho}(ii) and $\frac{1}{\infty}:=0$.
\end{theo}

\begin{theo} \label{pinnapplica}
 Assume {\bf (A)}. Let $u \in H^{2, \frac{d}{2}}(U) \cap H^{1,2}_0(U) \cap L^{\infty}(U)$ be a unique strong solution to \eqref{maineqpa} as in Theorem \ref{mainexisthm}.  Let $\Psi \in H^{2,2}(U) \cap H^{1,2}_0(U) \cap L^{\infty}(U)$, $f \in L^q(U) \cap L^{2}(U)$ and $c \in L^{\frac{d}{2}}(U) \cap L^2(U)$. Then, there exists 
$\mu=\rho dx$ with $\rho \in H^{1,p}(U) \cap C(\overline{U})$ and $\rho(x)>0$ for all $x \in \overline{U}$ such that
\begin{equation} \label{pinnesti}
\Vert u-\Psi \Vert_{L^2(U, \mu)} \leq \frac{1}{\alpha} \| f- \mathcal{L} \Psi \|_{L^2(U, \mu)},
\end{equation}
where 
\begin{equation} \label{lpsidefined}
\mathcal{L}\Psi = -\Delta \Psi + \langle \mathbf{H}, \nabla \Psi \rangle + (c+\alpha) \Psi.
\end{equation}
In particular,
\begin{equation} \label{pracpinestim}
\Vert u-\Psi \Vert_{L^2(U)}\leq \frac{\sqrt{K_1}}{\alpha} \| f- \mathcal{L} \Psi \|_{L^2(U)},
\end{equation}
where $K_1 \geq 1$ is a constant as in Theorem \ref{existrho}(ii).
\end{theo}
The main novelty of this paper is that, by identifying the solution as the image of a weighted resolvent associated with a Dirichlet form on $L^2(U, \rho\,dx)$, we obtain contraction estimates with improved and explicit constants. In \cite{L25jm}, since the coefficient $c$ is assumed to belong to $L^1(U)$, it is pointed out that the associated bilinear form need not be well-defined on all of $H_0^{1,2}(U)$. Focusing on the contraction estimates in \eqref{contrawithrele} of Theorem \ref{mainexisthm} and comparing them with \eqref{eltwoesljm}, one can observe that the constant improves from $K_1$ to $K_1^{\frac{1}{\theta}}$. In fact, the constant $K_1$ arises from the Harnack inequality and is a highly complicated constant. As a result, a careful and detailed calculation may yield a value of $K_1$ that is significantly larger than one. In particular, when $K_1$ is large and $\theta \ge 2$, the quantity $K_1^{\frac{1}{\theta}}$ is substantially smaller than $K_1$.
Moreover, in the case $\theta = \infty$, \eqref{contrawithrele} in Theorem \ref{mainexisthm} partially improves upon \cite[Chapter 11, Section 2, Theorem 1]{Kr08} for the case of bounded coefficients.
On the other hand, Theorem \ref{pinnapplica} concerns error estimates arising from the approximation of solutions to partial differential equations using artificial intelligence (AI) techniques. This result also partially improves the estimate in \cite[Remark 5.3]{L25jm},
\begin{equation*} 
\Vert u-\Psi \Vert_{L^2(U)} \leq \frac{K_1}{\alpha} \, \| f- \mathcal{L} \Psi \|_{L^2(U)},
\end{equation*}
since $\sqrt{K_1}$ is significantly smaller than $K_1$ when $K_1 \ge 1$ is sufficiently large. Further discussion of this point is provided in Remark \ref{finalrema}.\\ 
The structure of the paper is as follows. In Section \ref{sec2}, we present the notations and definitions used throughout the paper. Section \ref{sec3} constructs the Dirichlet form and the resolvent associated with \eqref{maineqpa}. Section \ref{sec4} contains the proofs of the main results.

\section{Notations and definitions} \label{sec2}
For $a,b \in \mathbb{R}$ with $a \geq b$, we define $a \vee b = b \vee a := a$ and $a \wedge b = b \wedge a := b$. We consider the Euclidean space $\mathbb{R}^d$ with $d \geq 1$, equipped with the standard inner product $\langle \cdot, \cdot \rangle$ and the corresponding Euclidean norm $\|\cdot\|$.  
For $x_0 \in \mathbb{R}^d$ and $r > 0$, define $B_r(x_0) := \{ x \in \mathbb{R}^d : \|x - x_0\| < r \}$. 
Let $U \subset \mathbb{R}^d$ be an open set. Denote the Lebesgue measure on $\mathbb{R}^d$ by $dx$. For a Lebesgue measurable set $V$, write $|V|:=dx(V)$. The spaces of continuous functions on $U$ and on its closure $\overline{U}$ are denoted by $C(U)$ and $C(\overline{U})$, respectively.  
We write $C_0(U)$ for the set of all continuous functions whose support is compact and contained in $U$.  
For $k \in \mathbb{N} \cup \{\infty\}$, the space of $k$-times continuously differentiable functions on $U$ is denoted by $C^k(U)$, and we define $C_0^k(U) := C^k(U) \cap C_0(U)$.
Let $s \in [1,\infty]$. The Lebesgue space with respect to a measure $\nu$, equipped with the standard $L^s$-norm, is denoted by $L^s(U,\nu)$.  
The space of $L^s$-integrable vector fields on $U$ with respect to $\nu$ is denoted by $L^s(U,\mathbb{R}^d,\nu)$, with norm
$\|F\|_{L^s(U,\nu)} := \big\| \|F \| \big \|_{L^s(U, \nu)}$, $F \in L^s(U,\mathbb{R}^d,\nu)$.
Define $L^s(U):=L^s(U, dx)$ and $L^s(U, \mathbb{R}^d):=L^s(U, \mathbb{R}^d, dx)$. For a function $f$, the weak derivative with respect to the $i$-th coordinate is denoted by $\partial_i f$.  
The Sobolev space $H^{1,s}(U)$ is the space of all $f \in L^s(U)$ such that $\partial_i f \in L^s(U)$ for each $i=1,\dots,d$ equipped with the usual $H^{1,s}(U)$-norm. We denote by $H^{1,2}_0(U)$ the closure of $C_0^{\infty}(U)$ with respect to the norm $\|\cdot\|_{H^{1,2}(U)}$.  
The Sobolev space $H^{2,s}(U)$ consists of all $f \in L^s(U)$ such that all first- and second-order weak derivatives belong to $L^s(U)$.  
Finally, for any twice differentiable function $f$, the Laplacian is defined as
$\Delta f := \sum_{i=1}^d \partial_i^2 f $.
 
\begin{defn} \label{defnsoluti}
Let $U$ be a bounded open subset of $\mathbb{R}^d$. Let $\rho \in L^{\infty}(U)$, $\mathbf{H} \in L^2(U, \mathbb{R}^d)$, $\mathbf{B} \in L^2(U, \mathbb{R}^d)$, $c \in L^1(U)$, $f \in L^1(U)$, and let $\alpha>0$ be a constant.  
Assume that $u \in H_0^{1,2}(U)$ and $cu \in L^1(U)$.  
\begin{itemize}
\item[(i)]  
We say that $u$ is a weak solution to \eqref{maineqpa} if the following identity holds:
\begin{equation} \label{maineqfir}
\int_U \langle \nabla u, \nabla \psi \rangle 
+ \langle \mathbf{H}, \nabla u \rangle \psi 
+ (c+\alpha) u \psi \, dx 
= \int_U f \psi \, dx
\quad \text{ for all } \psi \in C_0^\infty(U).
\end{equation}

\item[(ii)]  
We say that $u$ is a weak solution to \eqref{maineq3} if the following identity holds:
\begin{equation} \label{maineqfir2}
\int_U \langle \rho \nabla u, \nabla \varphi \rangle 
+ \langle \rho \mathbf{B}, \nabla u \rangle \varphi 
+ \rho (c+\alpha) u \varphi \, dx 
= \int_U \rho f \varphi \, dx
\quad \text{ for all } \varphi \in C_0^\infty(U).
\end{equation}

\item[(iii)]  
We say that $u$ is a strong solution to \eqref{maineqpa} if  
$u \in H^{2,1}(U) \cap H^{1,2}_0(U)$ and satisfies
\[
- \Delta u(x) 
+ \langle \mathbf{H}(x), \nabla u(x) \rangle 
+ \big(c(x)+\alpha\big) u(x) 
= f(x) \quad \text{for a.e. } x \in U.
\]
\end{itemize}
\end{defn}
\centerline{}
\noindent
\begin{defn}[Dirichlet form] \label{secdiri}
Let $\nu$ be a measure on $\mathbb{R}^d$.
A pair $(\mathcal E, D(\mathcal E))$ is called a Dirichlet form
on $L^2(\mathbb{R}^d,\nu)$ if the following conditions are satisfied:
\begin{itemize}
\item[(i)]
$D(\mathcal E)$ is a dense linear subspace of $L^2(\mathbb{R}^d,\nu)$.

\item[(ii)]
$\mathcal E : D(\mathcal E) \times D(\mathcal E) \to \mathbb{R}$ is a positive definite bilinear form.

\item[(iii)]
$D(\mathcal E)$ is complete with respect to the norm
\[
\|u\|_{\mathcal E_1}
:= \bigl( \mathcal E(u,u) + \|u\|_{L^2(\mathbb{R}^d,\nu)}^2 \bigr)^{1/2}.
\]

\item[(iv)]
$(\mathcal E, D(\mathcal E))$ satisfies the strong sector condition, i.e.,
there exists a constant $K \ge 1$ such that
\[
|\mathcal E(u,v)|
\le K\,\mathcal E(u,u)^{1/2}\,\mathcal E(v,v)^{1/2}
\quad \text{for all } u,v \in D(\mathcal E).
\]

\item[(v)]
For all $u \in D(\mathcal E)$,
\[
u^+ \wedge 1 \in D(\mathcal E),
\]
and
\[
\mathcal E(u + u^+ \wedge 1,\; u - u^+ \wedge 1) \ge 0,
\qquad
\mathcal E(u - u^+ \wedge 1,\; u + u^+ \wedge 1) \ge 0.
\]
\end{itemize}
\end{defn}

\section{Dirichlet forms and resolvents} \label{sec3}
The following constitutes the core result of \cite{L25jm}, whose original idea goes back to \cite[Theorem 1]{BRS12} and \cite[Chapter 2]{BKRS15}.
\begin{theo} \label{existrho} 
Assume that {\bf (A)} holds. Let $B_r(x_0)$ be an open ball satisfying $\overline{U} \subset B_r(x_0)$. Extend $\mathbf{H} \in L^p(U, \mathbb{R}^d)$ to $L^p(\mathbb{R}^d, \mathbb{R}^d)$ by zero extension. Then, the following hold:
\begin{itemize}
\item[(i)]
Let $x_1 \in U$. Then, there exists $\rho \in H^{1,2}(B_{4r}(x_0)) \cap C(B_{4r}(x_0))$ with $\rho(x)>0$ for all $x \in B_{4r}(x_0)$ and $\rho(x_1)=1$ such that
\begin{equation*} 
\int_{B_{4r}(x_0)} \langle   \nabla \rho  + \rho \mathbf{H}, \nabla \varphi \rangle dx = 0 \quad \text{ for all $\varphi \in C_0^{\infty}(B_{4r}(x_0))$}.
\end{equation*}
\item[(ii)]
Let $\rho$ be as in Theorem \ref{existrho}(i). Then, there exists a constant $K_1 \geq 1$ which only depends on $d$, $r$, $p$, $\|\mathbf{H}\|_{L^p(U)}$ such that \;$\max_{\overline{B}_{3r}(x_0)} \rho  \leq K_1 \min_{\overline{B}_{3r}(x_0)} \rho$, and hence
$$
1 \leq \max_{\overline{U}} \rho  \leq K_1 \min_{\overline{U}} \rho \leq K_1.
$$
\item[(iii)]
We have $\rho \in H^{1,p}(B_r(x_0))$.
\end{itemize}
\end{theo}
\begin{proof}
Assertions (i) and (ii) follow from \cite[Theorem 3.1(i), (ii)]{L25jm}. The assertion (iii) holds by \cite[Theorem 4.1]{L25jm}.
\end{proof} 

\begin{theo} \label{divfretransform}
Suppose that assumption {\bf (A)} is satisfied.  
Let $\rho \in H^{1,2}(U) \cap C(\overline{U})$ be a strictly positive function on $\overline{U}$, obtained as in Theorem \ref{existrho}.  
Define the vector field
\begin{equation} \label{divfreebdef}
\mathbf{B} := \mathbf{H} + \frac{1}{\rho} \nabla \rho \quad \text{in } U.
\end{equation}
Then $\rho \mathbf{B} \in L^2(U, \mathbb{R}^d)$ and it holds that
\begin{equation} \label{divfeepropbolb}
\int_U \langle \rho \mathbf{B}, \nabla \varphi \rangle \, dx = 0 
\quad \text{for every } \varphi \in C_0^\infty(U).
\end{equation}
Let $f \in L^1(U)$ and $u \in H^{1,2}_0(U)$ with $cu \in L^1(U)$. Then, $u$ satisfies \eqref{maineqfir}, if and only if $u$ satisfies \eqref{maineqfir2}. Equivalently, $u$ is a weak solution to \eqref{maineqpa} if and only if $u$ is a weak solution to \eqref{maineq3}.
\end{theo}
\begin{proof}
The argument follows the same lines as in \cite[Theorem 3.2]{L25jm}, 
with the special case $\mathbf{F} = 0$.  For the reader’s convenience, we include the
details below.\\
First, assume that $u$ is a weak solution to \eqref{maineqpa}. By a standard
density argument, $u$ satisfies
\begin{equation} \label{maineqfirexte}
\int_U \langle \nabla u, \nabla \psi \rangle
+ \langle \mathbf{H}, \nabla u \rangle \psi
+ (c+\alpha) u \psi \, dx
= \int_U f \psi \, dx
\quad \text{for all } \psi \in H^{1,2}_0(U) \cap L^{\infty}(U).
\end{equation}
Let $\varphi \in C_0^{\infty}(U)$ be arbitrary. By the product rule (see, e.g.,
\cite[Chapter I, Corollary 4.15]{MR92}), we have
$\rho \varphi \in H^{1,2}_0(U) \cap L^{\infty}(U)$. Choosing
$\psi = \rho \varphi$ in \eqref{maineqfirexte}, we obtain
\begin{align*}
\int_U f \varphi \rho \, dx
&= \int_U \langle \nabla u, \nabla(\rho \varphi) \rangle
   + \langle \mathbf{H}, \nabla u \rangle \rho \varphi
   + (c+\alpha) u \rho \varphi \, dx \\
&= \int_U \langle \rho \nabla u, \nabla \varphi \rangle
   + \Big\langle \frac{1}{\rho} \nabla \rho + \mathbf{H}, \nabla u \Big\rangle
     \rho \varphi
   + \rho (c+\alpha) u \varphi \, dx \\
&= \int_U \langle \rho \nabla u, \nabla \varphi \rangle
   + \langle \rho \mathbf{B}, \nabla u \rangle \varphi
   + \rho (c+\alpha) u \varphi \, dx .
\end{align*}
Hence, $u$ is a weak solution to \eqref{maineq3}.\\
Conversely, assume that $u$ is a weak solution to \eqref{maineq3}. By a density
argument, $u$ satisfies
\begin{equation} \label{maineqfir2exte}
\int_U \langle \rho \nabla u, \nabla \varphi \rangle
+ \langle \rho \mathbf{B}, \nabla u \rangle \varphi
+ \rho (c+\alpha) u \varphi \, dx
= \int_U \rho f \varphi \, dx
\quad \text{for all } \varphi \in H^{1,2}_0(U) \cap L^{\infty}(U).
\end{equation}
Let $\psi \in C_0^{\infty}(U)$ be arbitrary. By the product rule,
$\frac{\psi}{\rho} \in H^{1,2}_0(U) \cap L^{\infty}(U)$. Taking
$\varphi = \frac{\psi}{\rho}$ in \eqref{maineqfir2exte}, we obtain
\begin{align*}
\int_U f \psi \, dx
&= \int_U \Big\langle \rho \nabla u,
      \nabla\!\left(\frac{\psi}{\rho}\right) \Big\rangle
   + \langle \mathbf{B}, \nabla u \rangle \psi
   + (c+\alpha) u \psi \, dx \\
&= \int_U \langle \nabla u, \nabla \psi \rangle
   + \Big\langle -\frac{1}{\rho}  \nabla \rho + \mathbf{B}, \nabla u \Big\rangle
     \psi
   + (c+\alpha) u \psi \, dx \\
&= \int_U \langle  \nabla u, \nabla \psi \rangle
   + \langle \mathbf{H}, \nabla u \rangle \psi
   + (c+\alpha) u \psi \, dx .
\end{align*}
Therefore, $u$ is a weak solution to \eqref{maineqpa}.

\end{proof}
\centerline{}
Under assumption {\bf (A)}, let $\rho$ be the function provided by Theorem \ref{existrho}, and set
\begin{align} \label{rhodivfreeb}
\mu := \rho \, dx.
\end{align}
Let $\mathbf{B}$ be defined as in \eqref{divfreebdef}. From Theorem \ref{existrho}, we have $\rho \mathbf{B} \in L^p(U, \mathbb{R}^d)$ and
\begin{equation*} 
\int_{U} \langle \mathbf{B}, \nabla \varphi \rangle \, d\mu
= \int_{U} \langle \rho \mathbf{B}, \nabla \varphi \rangle \, dx
= \int_{U} \langle \rho \mathbf{H} +  \nabla \rho, \nabla \varphi \rangle \, dx 
= 0\;
\quad \text{ for all } \, \varphi \in C_0^{\infty}(U),
\end{equation*}
so that \eqref{divfeepropbolb} follows. In particular, by the product rule and \eqref{divfeepropbolb},
\begin{equation} \label{divfreeconalg}
\int_{U} \langle \mathbf{B}, \nabla \varphi \rangle \, \varphi \, d\mu
= \frac12 \int_{U} \langle \mathbf{B}, \nabla (\varphi^2) \rangle \, d\mu
= 0
\quad \text{  for all } \, \varphi \in C_0^{\infty}(U).
\end{equation}
Moreover, by the product rule and \eqref{divfeepropbolb}, we have
\begin{equation} \label{skewsympro}
\int_{U} \langle \mathbf{B}, \nabla \varphi \rangle \psi\, d\mu= - \int_{U} \langle \mathbf{B}, \nabla \psi \rangle \varphi \, d\mu \quad \text{ for all $\varphi, \psi \in C_0^{\infty}(U)$}.
\end{equation}
Let $f, g \in C_0^{\infty}(U)$.  
Since $\rho$ is bounded above and below by strictly positive constants on $\overline{U}$ by Theorem \ref{existrho}(ii), the H\"{o}lder inequality and the Sobolev inequality (see \cite[Section 5.6, Theorem 1]{Ev10})
yield
\begin{align}
\left| \int_{U} \langle \mathbf{B}, \nabla f \rangle g \, d\mu \right|
&\leq \|\rho \mathbf{B}\|_{L^{d}(U)} \|\nabla f\|_{L^2(U)} \|g\|_{L^{\frac{2d}{d-2}}(U)} \nonumber \\
&\leq \|\rho \mathbf{B}\|_{L^{d}(U)} \|\nabla f\|_{L^2(U)} \cdot \frac{2(d-1)}{d-2}  \|\nabla g\|_{L^{2}(U)} \nonumber \\
&\leq  \frac{\max_{\overline{U}} \rho}{\min_{\overline{U}} \rho} \|\mathbf{B}\|_{L^d(U)} \frac{2(d-1)}{d-2}
\| \nabla f \|_{L^2(U, \mu)}\| \nabla g \|_{L^2(U, \mu)}. \label{weksecb}
\end{align}
Moreover, by the H\"{o}lder inequality and the Sobolev inequality (see \cite[Section 5.6, Theorem 1]{Ev10}),
\begin{align}
\left |\int_{U}c  f g\, d\mu \right| &\leq  \max_{\overline{U}} \rho \cdot \|c \|_{L^{\frac{d}{2}}(U)} \|f\|_{L^{\frac{2d}{d-2}}(U)} \|g\|_{L^{\frac{2d}{d-2}}(U)} \nonumber \\
& \leq  \max_{\overline{U}} \rho \cdot  \|c\|_{L^{\frac{d}{2}}(U)} \left(  \frac{2(d-1)}{d-2}	\right)^2 \| \nabla f \|_{L^2(U)}\| \nabla g \|_{L^2(U)} \nonumber \\
& \leq \frac{\max_{\overline{U}} \rho}{\min_{\overline{U}} \rho} \cdot  \|c\|_{L^{\frac{d}{2}}(U)} \left(  \frac{2(d-1)}{d-2}	\right)^2 \| \nabla f \|_{L^2(U, \mu)}\| \nabla g \|_{L^2(U, \mu)}. \label{weksecbzer}
\end{align}
We now define the positive bilinear form $(\mathcal{E}, C_0^{\infty}(U))$ by
\begin{equation} \label{predirici}
\mathcal{E}(f,g)
:= \int_{U} \langle \nabla f, \nabla g \rangle \, d\mu
   + \int_{U} \langle \mathbf{B}, \nabla f \rangle g \, d\mu + \int_{U}c fg \,d\mu,
\quad f, g \in C_0^{\infty}(U).
\end{equation}
The next result shows that $(\mathcal{E}, C_0^{\infty}(U))$ is closable in $L^2(U, \mu)$ and that its closure is a Dirichlet form (see \cite[Chapter I]{MR92}). The basic structure of the proof is analogous to that of \cite[Theorem 3.2]{L25bv}.
\begin{prop} \label{matebisdiri}
Assume {\bf (A)} holds.  
Let $\rho \in H^{1,p}(U) \cap C(\overline{U})$ be as in Theorem \ref{existrho}, and define $\mathbf{B}$ and $\mu$ according to \eqref{divfreebdef} and \eqref{rhodivfreeb}, respectively.  
Consider the bilinear form $(\mathcal{E}, C_0^{\infty}(U))$ specified in \eqref{predirici}.  
Then the following hold:

\begin{itemize}
\item[(i)]  
The bilinear form $(\mathcal{E}, C_0^{\infty}(U))$ satisfies the strong sector condition, i.e. there exists a constant $c>0$, independent of $f,g \in C_0^{\infty}(U)$, such that  
\[
|\mathcal{E}(f,g)| \leq c\,\mathcal{E}(f,f)^{1/2}\,\mathcal{E}(g,g)^{1/2} \quad \text{ for all } f,g \in C_0^{\infty}(U).
\]

\item[(ii)]  
The bilinear form $(\mathcal{E}, C_0^{\infty}(U))$ is closable in $L^2(U,\mu)$ (see \cite[Chapter I, Definition 3.1]{MR92}). Moreover, its closure, denoted by $(\mathcal{E}, D(\mathcal{E}))$, is a coercive closed form on $L^2(U,\mu)$ (cf. \cite[Chapter I, Definition 2.4, Remark 3.2]{MR92}). 

\item[(iii)]  
$D(\mathcal{E}) = H^{1,2}_0(U)$ and
\[
\mathcal{E}(f,g) = \int_{U} \langle \nabla f, \nabla g \rangle \, d\mu
+ \int_{U} \langle \mathbf{B}, \nabla f \rangle\, g \, d\mu
+ \int_{U} c\, f g \, d\mu \quad \text{ for all } f,g \in D(\mathcal{E}).
\]

\item[(iv)]  
$(\mathcal{E}, D(\mathcal{E}))$ is a Dirichlet form (see Definition \ref{secdiri}).
\end{itemize}
\end{prop}
\begin{proof}
(i)
Let $f, g \in C_0^{\infty}(U)$.  
By applying the Cauchy-Schwarz inequality together with \eqref{weksecb} and \eqref{weksecbzer}, we obtain
\begin{align*}
\left| \mathcal{E}(f,g) \right| 
&\le \left| \int_{U} \langle \nabla f, \nabla g \rangle \, d\mu \right| 
     + \left| \int_{U} \langle \mathbf{B}, \nabla f \rangle \, g \, d\mu \right|
     + \left| \int_{U} c f g \, d\mu \right| \\
&\le \mathcal{E}(f,f)^{1/2} \, \mathcal{E}(g,g)^{1/2} 
     + \frac{\max_{\overline{U}} \rho}{\min_{\overline{U}} \rho} 
       \| \mathbf{B} \|_{L^d(U)} \gamma \, \mathcal{E}(f,f)^{1/2} \, \mathcal{E}(g,g)^{1/2} \\
&\qquad + \frac{\max_{\overline{U}} \rho}{\min_{\overline{U}} \rho} 
       \|c\|_{L^{\frac{d}{2}}(U)} \, \gamma^2 \, \mathcal{E}(f,f)^{1/2} \, \mathcal{E}(g,g)^{1/2} \\
&= \left( 1 +  \frac{\max_{\overline{U}} \rho}{\min_{\overline{U}} \rho} 
           \| \mathbf{B} \|_{L^d(U)}  \gamma
         + \frac{\max_{\overline{U}} \rho}{\min_{\overline{U}} \rho} 
           \|c\|_{L^{\frac{d}{2}}(U)} \, \gamma^2 \right) 
   \mathcal{E}(f,f)^{1/2} \, \mathcal{E}(g,g)^{1/2},
\end{align*}
where $\gamma = \frac{2(d-1)}{d-2}$. Thus, the assertion follows. \\ \\
(ii) To prove that $(\mathcal{E}, C_0^{\infty}(U))$ is closable, suppose $(f_n)_{n \ge 1} \subset C_0^{\infty}(U)$ satisfies  
$\lim_{n \to \infty} f_n = 0$ in $L^2(U, \mu)$ and  
$\lim_{n,m \to \infty} \mathcal{E}(f_n - f_m, f_n - f_m) = 0$.  
From \eqref{divfreeconalg}, we have
\begin{align}
&\int_{U} \|\nabla f_n - \nabla f_m\|^2 \, dx \leq \frac{1}{\min_{\overline{U}} \rho} \int_{U} \|\nabla f_n - \nabla f_m\|^2 \, d\mu \nonumber  \\
&\le \frac{1}{\min_{\overline{U}} \rho}\left( \int_{U} \|\nabla(f_n - f_m)\|^2 \, d\mu 
     + \frac12 \int_{U} \langle \mathbf{B}, \nabla (f_n - f_m)^2 \rangle \, d\mu 
     + \int_{U} c (f_n - f_m)^2 \, d\mu \right)  \nonumber \\
&= \frac{1}{\min_{\overline{U}} \rho} \mathcal{E}(f_n - f_m, f_n - f_m) \longrightarrow 0
\quad \text{as } n, m \to \infty. \label{diribdsobno}
\end{align}
Thus,
\[
\|\nabla(f_n - f_m)\|_{L^2(U)}^2 \longrightarrow 0 \quad \text{as } n,m \to \infty.
\]
Since $L^2(U)$ is complete, for each $i \in \{1, \dots, d\}$ there exists $g_i \in L^2(U)$ such that
\[
\partial_i f_n \longrightarrow g_i \quad \text{in } L^2(U) \quad \text{as } n \to \infty.
\]
Now, let $\varphi \in C_0^{\infty}(U)$ and fix $i \in \{1, \dots, d\}$. Then
\[
\int_{U} g_i \varphi \, dx
= \lim_{n \to \infty} \int_{U} \partial_i f_n \, \varphi \, dx
= - \lim_{n \to \infty} \int_{U} f_n \, \partial_i \varphi \, dx
= 0.
\]
Therefore $g_i = 0$ in $U$, implying $\nabla f_n \to 0$ in $L^2(U,\mathbb{R}^d)$.
Using \eqref{divfreeconalg} together with the H\"{o}lder inequality and the Sobolev inequality, we deduce that
\begin{align*}
\mathcal{E}(f_n, f_n) 
&= \int_{U} \|\nabla f_n \|^2 \, d\mu
   + \frac12 \int_{U} \langle \mathbf{B}, \nabla (f_n^2) \rangle \, d\mu + \int_{U} c (f_n)^2 \, d\mu  \\
&\leq \max_{\overline{U}} \rho \left( \int_{U} \|\nabla f_n \|^2 \, dx
   + \int_{U} c (f_n)^2 \, dx \right) \\
&\leq \max_{\overline{U}} \rho \left( 1 + \|c\|_{L^{\frac{d}{2}}(U)} 
   \left( \frac{2(d-1)}{d-2} \right)^2 \right)
   \| \nabla f_n  \|_{L^2(U)}^2,
\end{align*}
and hence $\mathcal{E}(f_n,f_n) \rightarrow 0$ as $n \to \infty$.  
This shows that $(\mathcal{E}, C_0^{\infty}(U))$ is closable.  
The remaining assertions follow directly from the sector condition established in (i). \\ \\
(iii) 
Let $f \in D(\mathcal{E})$.  
Then there exists a sequence $(f_n)_{n \ge 1} \subset C_0^{\infty}(U)$ such that  
$f_n \to f$ in $L^2(U, \mu)$ and  
$\mathcal{E}(f_n - f_m, f_n - f_m) \to 0$ as $n, m \to \infty$.  
By \eqref{diribdsobno}, we have  
$\| \nabla (f_n - f_m) \|_{L^{2}(U)} \to 0$ as $n, m \to \infty$.  
Since $H^{1,2}_0(U)$ is complete, there exists $\tilde{f} \in H^{1,2}_0(U)$ such that  
$f_n \to \tilde{f}$ in $H^{1,2}_0(U)$, and thus $f = \tilde{f} \in H^{1,2}_0(U)$.
Conversely, assume $f \in H^{1,2}_0(U)$.  
Then there exists $(f_n)_{n \ge 1} \subset C_0^{\infty}(U)$ with $f_n \to f$ in $H^{1,2}_0(U)$,  
which in particular implies $f_n \to f$ in $L^2(U, \mu)$.  
Using \eqref{divfreeconalg} together with the H\"{o}lder inequality and the Sobolev inequality, we obtain
\begin{align*}
\mathcal{E}(f_n - f_m, f_n - f_m) 
&= \int_{U} \|\nabla(f_n - f_m)\|^2 \, d\mu
   + \frac12 \int_{U} \langle \mathbf{B}, \nabla (f_n - f_m)^2 \rangle \, d\mu + \int_{U} c (f_n - f_m)^2 \, d\mu  \\
&\leq  \max_{\overline{U}} \rho \left( \int_{U} \|\nabla(f_n - f_m)\|^2 \, dx
   + \int_{U} c (f_n - f_m)^2 \, dx \right) \\
&\le \max_{\overline{U}} \rho \left( 1 + \|c\|_{L^{\frac{d}{2}}(U)} 
   \left( \frac{2(d-1)}{d-2} \right)^2 \right)
   \| \nabla(f_n - f_m) \|_{L^2(U)}^2,
\end{align*}
which tends to $0$ as $n, m \to \infty$.  
Hence $f \in D(\mathcal{E})$. For the remaining statement, take $f, g \in D(\mathcal{E})$.  
There exist sequences $(f_n)_{n \ge 1}, (g_n)_{n \ge 1} \subset C_0^{\infty}(U)$ such that  
$f_n \to f$ and $g_n \to g$ in $L^2(U, \mu)$, and  
$\mathcal{E}(f_n - f_m, f_n - f_m) \to 0$ and  
$\mathcal{E}(g_n - g_m, g_n - g_m) \to 0$ as $n, m \to \infty$.  
From the estimate in \eqref{diribdsobno} and the completeness of $H^{1,2}_0(U)$, it follows that  
$f, g \in H^{1,2}_0(U)$ and $f_n \to f$, $g_n \to g$ in $H^{1,2}_0(U)$.  
Finally, since $\mathcal{E}(f, g) = \lim_{n \to \infty} \mathcal{E}(f_n, g_n)$,  
the claim is proved. \\ \\
(iv) Let $\varepsilon > 0$ and choose $\varphi_{\varepsilon} \in C^{\infty}(\mathbb{R})$ as in \cite[Proposition 5.9]{L25bv}.  
First, take $u \in C_0^{\infty}(U)$. Then $\varphi_{\varepsilon} \circ u \in C_0^{\infty}(U)$.  
Consider the bilinear form $(\mathcal{E}^{0}, C_0^{\infty}(U))$ given by
\begin{equation} \label{secdirich}
\mathcal{E}^{0}(f,g) := \int_{U} \langle \nabla f, \nabla g \rangle \, d\mu,
\quad f,g \in C_0^{\infty}(U).
\end{equation}
A direct computation shows that
\[
\mathcal{E}^{0}(\varphi_{\varepsilon} \circ u,\, u - \varphi_{\varepsilon} \circ u)
= \int_{U} (\varphi'_{\varepsilon} \circ u)\,(1 - \varphi'_{\varepsilon} \circ u) \, \|\nabla u\|^2 \, d\mu \ge 0,
\]
which implies $\liminf_{\varepsilon \to 0+} \mathcal{E}^{0}(\varphi_{\varepsilon} \circ u,\, u - \varphi_{\varepsilon} \circ u) \ge 0$.  
Since 
$$
\mathcal{E}^{0}(u - \varphi_{\varepsilon} \circ u,\, \varphi_{\varepsilon} \circ u) = \mathcal{E}^{0}(\varphi_{\varepsilon} \circ u,\, u - \varphi_{\varepsilon} \circ u), 
$$
it follows that
\[
\liminf_{\varepsilon \to 0+} \mathcal{E}^{0}(u - \varphi_{\varepsilon} \circ u,\, \varphi_{\varepsilon} \circ u) \ge 0.
\]
Now define
\[
S_{\varepsilon} := \int_{U} \langle \mathbf{B},\, \nabla (\varphi_{\varepsilon} \circ u) \rangle \,(u - \varphi_{\varepsilon} \circ u) \, d\mu.
\]
By the Lebesgue dominated convergence theorem and \cite[Proposition 5.10]{L25bv}, we obtain
\begin{align}
\lim_{\varepsilon \to 0+} S_{\varepsilon}
&= \lim_{\varepsilon \to 0+} \int_{U} (\varphi'_{\varepsilon} \circ u)\,(u - \varphi_{\varepsilon} \circ u) \, \langle \mathbf{B}, \nabla u \rangle \, d\mu \nonumber \\
&= \int_{U} (1_{[0,1]} \circ u)\,(u - u^+ \wedge 1) \, \langle \mathbf{B}, \nabla u \rangle \, d\mu = 0. \label{svarepsilon}
\end{align}
Similarly, let
\[
T_{\varepsilon} := \int_{U} \langle \mathbf{B},\, \nabla (u - \varphi_{\varepsilon} \circ u) \rangle \, (\varphi_{\varepsilon} \circ u) \, d\mu.
\]
By \eqref{skewsympro}, we have $T_{\varepsilon} = -S_{\varepsilon}$, and thus by \eqref{svarepsilon},
\begin{equation} \label{tvarepsilon}
\lim_{\varepsilon \to 0+} T_{\varepsilon} = 0.
\end{equation}
Combining \eqref{secdirich}, \eqref{svarepsilon} and \eqref{tvarepsilon}, we find
\begin{align*}
&\liminf_{\varepsilon \to 0+} \mathcal{E}(\varphi_{\varepsilon} \circ u,\, u - \varphi_{\varepsilon} \circ u) \\
&= \liminf_{\varepsilon \to 0+} \left( \mathcal{E}^{0}(\varphi_{\varepsilon} \circ u,\, u - \varphi_{\varepsilon} \circ u) + S_{\varepsilon}  + \int_{U} c (\varphi_{\varepsilon} \circ u) (u - \varphi_{\varepsilon} \circ u )\,d\mu\right) \\
&\ge \liminf_{\varepsilon \to 0+} \mathcal{E}^{0}(\varphi_{\varepsilon} \circ u,\, u - \varphi_{\varepsilon} \circ u)
   + \liminf_{\varepsilon \to 0+} S_{\varepsilon}  + \liminf_{\varepsilon \rightarrow 0+}\int_{U} c (\varphi_{\varepsilon} \circ u) (u - \varphi_{\varepsilon} \circ u )\,d\mu\\
&= \liminf_{\varepsilon \to 0+} \mathcal{E}^{0}(\varphi_{\varepsilon} \circ u,\, u - \varphi_{\varepsilon} \circ u)  + \int_{U} c (u^+\wedge1)(u-u^+ \wedge 1)\,d\mu \ge 0.
\end{align*}
Likewise, using \eqref{tvarepsilon},
\begin{align*}
&\liminf_{\varepsilon \to 0+} \mathcal{E}(u - \varphi_{\varepsilon} \circ u,\, \varphi_{\varepsilon} \circ u) \\
&= \liminf_{\varepsilon \to 0+} \left( \mathcal{E}^{0}(u - \varphi_{\varepsilon} \circ u,\, \varphi_{\varepsilon} \circ u) + T_{\varepsilon} +\int_{U} c (\varphi_{\varepsilon} \circ u) (u - \varphi_{\varepsilon} \circ u )\,d\mu\right) \\
&\ge \liminf_{\varepsilon \to 0+} \mathcal{E}^{0}(u - \varphi_{\varepsilon} \circ u,\, \varphi_{\varepsilon} \circ u)
   + \liminf_{\varepsilon \to 0+} T_{\varepsilon} + \liminf_{\varepsilon \rightarrow 0+} \int_{U} c (\varphi_{\varepsilon} \circ u) (u - \varphi_{\varepsilon} \circ u )\,d\mu \\
&= \liminf_{\varepsilon \to 0+} \mathcal{E}^{0}(u - \varphi_{\varepsilon} \circ u,\, \varphi_{\varepsilon} \circ u) + \int_{U} c (u^+\wedge1)(u-u^+ \wedge 1)\,d\mu \ge 0.
\end{align*}
Since $C_0^{\infty}(U)$ is dense in $D(\mathcal{E})$ with respect to the norm 
$\mathcal{E}(\cdot,\cdot)^{1/2} + \|\cdot\|_{L^2(U,\mu)}$,  
it follows from \cite[Chapter I, Proposition 4.10]{MR92} that $(\mathcal{E}, D(\mathcal{E}))$ is a Dirichlet form.
\end{proof}
\centerline{}
Similarly to \eqref{predirici}, define the positive bilinear form $(\hat{\mathcal{E}}, C_0^{\infty}(U))$ by
\begin{equation} \label{prediricidua}
\hat{\mathcal{E}}(f,g)
:= \int_{U} \langle \nabla f, \nabla g \rangle \, d\mu
    +\int_{U} \langle -\mathbf{B}, \nabla f \rangle g \, d\mu + \int_{U}c fg \,d\mu,
\quad f, g \in C_0^{\infty}(U).
\end{equation}
Analogously to the proof of Proposition \ref{matebisdiri}, we obtain that $(\hat{\mathcal{E}}, C_0^{\infty}(U))$ is closable in $L^2(U, \mu)$ and that its closure is also a Dirichlet form. 

\begin{prop} \label{matebisdiridua}
Assume that {\bf (A)} holds.  
Let $\rho \in H^{1,p}(U) \cap C(\overline{U})$ be as in Theorem \ref{existrho}, and define $\mathbf{B}$ and $\mu$ by \eqref{divfreebdef} and \eqref{rhodivfreeb}, respectively.  
Consider the bilinear form $(\hat{\mathcal{E}}, C_0^{\infty}(U))$ given in \eqref{prediricidua}.  
Then the following statements hold:

\begin{itemize}
\item[(i)]  
The bilinear form $(\hat{\mathcal{E}}, C_0^{\infty}(U))$ satisfies the strong sector condition, i.e. there exists a constant $c>0$, independent of $f,g \in C_0^{\infty}(U)$, such that  
\[
|\hat{\mathcal{E}}(f,g)| \leq c\, \hat{\mathcal{E}}(f,f)^{1/2}\, \hat{\mathcal{E}}(g,g)^{1/2} 
\quad \text{ for all } f,g \in C_0^{\infty}(U).
\]

\item[(ii)]  
The bilinear form $(\hat{\mathcal{E}}, C_0^{\infty}(U))$ is closable in $L^2(U,\mu)$ (see \cite[Chapter I, Definition 3.1]{MR92}). Moreover, its closure, denoted by $(\hat{\mathcal{E}}, D(\hat{\mathcal{E}}))$, is a coercive closed form on $L^2(U,\mu)$ (cf. \cite[Chapter I, Definition 2.4, Remark 3.2]{MR92}). 

\item[(iii)]  
$D(\hat{\mathcal{E}}) = H^{1,2}_0(U)$ and
\[
\hat{\mathcal{E}}(f,g) 
= \int_{U} \langle \nabla f, \nabla g \rangle \, d\mu
+ \int_{U} \langle -\mathbf{B}, \nabla f \rangle\, g \, d\mu
+ \int_{U} c\, f g \, d\mu
\quad \text{ for all } f,g \in D(\mathcal{E}).
\]

\item[(iv)]  
$(\hat{\mathcal{E}}, D(\hat{\mathcal{E}}))$ is a Dirichlet form (see \cite[Chapter I, Definition 4.5]{MR92}).

\item[(v)]
For any $f,g \in D(\mathcal{E}) = D(\hat{\mathcal{E}}) = H^{1,2}_0(U)$,
\[
\hat{\mathcal{E}}(f,g) = \mathcal{E}(g,f),
\]
where $(\mathcal{E}, D(\mathcal{E}))$ is a Dirichlet form as in Proposition \ref{matebisdiri}.
\end{itemize}
\end{prop}
\begin{proof}
The proof of (i)-(iv) proceeds analogously to that of Proposition \ref{matebisdiri}(i)-(iv), with $\mathbf{B}$ in Proposition \ref{matebisdiri} replaced by $-\mathbf{B}$ here.  
We now prove (v). Let $f,g \in H^{1,2}_0(U)$.  
By the Sobolev inequality and the product rule, we have $fg \in H^{1,\frac{d}{d-1}}_0(U)$ and
$$
\nabla (fg) = g \nabla f + f \nabla g \quad \text{in } L^{\frac{d}{d-1}}(U, \mathbb{R}^d).
$$
It then follows from \eqref{divfeepropbolb} that
$$
\int_{U} \langle \mathbf{B}, \nabla(fg) \rangle \, d\mu = 0,
$$
and hence
\begin{equation*} 
\int_{U} \langle \mathbf{B}, \nabla g \rangle f \, d\mu
= - \int_{U} \langle \mathbf{B}, \nabla f \rangle g \, d\mu.
\end{equation*}
Therefore, $\hat{\mathcal{E}}(f,g) = \mathcal{E}(g,f)$, as desired.
\end{proof}

\begin{prop} \label{resolvpropo}
Let $(G_{\alpha})_{\alpha>0}$ be the strongly continuous sub-Markovian resolvent of contractions on $L^2(U, \mu)$ associated with the Dirichlet form $(\mathcal{E}, D(\mathcal{E}))$ as in Proposition \ref{matebisdiri}.
Then, for any $h \in L^{\infty}(U)$, $s \in [1, \infty]$, and $\alpha>0$, we have
\begin{equation} \label{elscontracpro}
\|G_{\alpha} h \|_{L^s(U, \mu)} \leq \frac{1}{\alpha} \|h\|_{L^s(U, \mu)}.
\end{equation}
\end{prop}
\begin{proof}
Let $\alpha>0$ be arbitrary.  
First, by \cite[Chapter I, Theorem 2.8]{MR92},  
\begin{equation} \label{firstresolid1}
\mathcal{E}(G_{\alpha}f, u) + \alpha \int_{U} G_{\alpha} f \, u \,d\mu
= \int_{U} f\, u \, d\mu 
\quad \text{ for all } f \in L^2(U, \mu), \ u \in D(\mathcal{E}).
\end{equation}
Let $(\hat{G}_{\alpha})_{\alpha>0}$ be the strongly continuous sub-Markovian resolvent of contractions on $L^2(U, \mu)$ associated with the Dirichlet form $(\hat{\mathcal{E}}, D(\hat{\mathcal{E}}))$ as in Proposition \ref{matebisdiridua}. Then, again by \cite[Chapter I, Theorem 2.8]{MR92},
\begin{equation} \label{firstresolid2}
\hat{\mathcal{E}}(\hat{G}_{\alpha}g, u) + \alpha \int_{U} \hat{G}_{\alpha} g \, u \,d\mu
= \int_{U} g\, u \, d\mu 
\quad \text{ for all } g \in L^2(U, \mu), \ u \in D(\hat{\mathcal{E}}).
\end{equation}
Now let $f,g \in L^{2}(U)$ and $\alpha>0$.  
Since $G_{\alpha}f, \hat{G}_{\alpha}g \in D(\mathcal{E}) = D(\hat{\mathcal{E}}) = H^{1,2}_0(U)$, combining \eqref{firstresolid1}, \eqref{firstresolid2}, and Proposition \ref{matebisdiridua}(v) yields
\begin{align}
\int_{U} f\, \hat{G}_{\alpha}g\,d\mu
&= \mathcal{E}(G_{\alpha}f, \hat{G}_{\alpha}g) 
+ \alpha \int_{U} G_{\alpha} f \, \hat{G}_{\alpha}g \,d\mu \nonumber \\
&= \hat{\mathcal{E}}(\hat{G}_{\alpha}g, G_{\alpha}f) 
+ \alpha \int_{U} \hat{G}_{\alpha}g \, G_{\alpha} f \, d\mu \nonumber \\
&= \int_{U} g\, G_{\alpha} f \, d\mu. \label{dualresolpro}
\end{align}
By \eqref{dualresolpro} and the sub-Markovian property of $(\hat{G}_{\alpha})_{\alpha>0}$, we have
\[
\int_{U} |G_{\alpha}f|\,d\mu
\leq \int_{U} G_{\alpha} |f| \,d\mu
= \int_{U} G_{\alpha} |f| \cdot 1_U \,d\mu
= \int_{U} |f| \cdot \hat{G}_{\alpha} 1_U\,d\mu
\leq \frac{1}{\alpha} \int_{U} |f| \,d\mu.
\]
By a standard argument using Cauchy sequences, $G_{\alpha}|_{L^2(U, \mu)}$ extends to an operator $G^{(1)}_{\alpha}$ satisfying
\begin{equation} \label{elonecontrac}
\|G^{(1)}_{\alpha} h\|_{L^1(U, \mu)} \leq \frac{1}{\alpha} \|h\|_{L^1(U, \mu)}
\quad \text{ for all } h \in L^1(U, \mu).
\end{equation}
Moreover, by the sub-Markovian property of $(G_{\alpha})_{\alpha>0}$
\begin{equation} \label{submarkres}
\|G^{(1)}_{\alpha} h\|_{L^{\infty}(U, \mu)}
\leq \frac{1}{\alpha} \|h\|_{L^{\infty}(U, \mu)}
\quad \text{ for all } h \in L^{\infty}(U, \mu).
\end{equation}
Combining \eqref{elonecontrac} with \eqref{submarkres} and applying the Riesz–Thorin interpolation theorem (\cite[Chapter 2, Theorem 2.1]{S11}), we obtain, for each $s \in [1, \infty]$,
\[
\|G^{(1)}_{\alpha} h\|_{L^{s}(U, \mu)} \leq \frac{1}{\alpha} \|h\|_{L^{s}(U, \mu)}
\quad \text{ for all } h \in L^{s}(U, \mu).
\]
Since $G_{\alpha} h = G_{\alpha}^{(1)}h$ for all $h \in L^{\infty}(U)$, \eqref{elscontracpro} follows.
\end{proof}
\centerline{}

\section{Proofs of main results} \label{sec4}
\noindent
This section provides the proofs of main results, Theorem \ref{mainexisthm} and Theorem \ref{pinnapplica}. 
\centerline{}
\centerline{}
{\bf Proof of Theorem \ref{mainexisthm}}\\
Fix a constant $\alpha > 0$. By \cite[Theorem 1.2]{L25jm}, there exists $u \in H^{2, \frac{d}{2}}(U) \cap H^{1,2}_0(U) \cap L^{\infty}(U)$ such that $u$ is a unique strong solution to \eqref{maineqpa}. Let $\theta \in [1, \infty]$. For each $n \geq 1$, define
\[
f_n := (f \wedge n) \vee (-n) \in L^{\infty}(U).
\]
For $\theta<\infty$, the Lebesgue dominated convergence theorem yields $f_n\to f$ in $L^{q\vee\theta}(U)$. For $\theta=\infty$, we have $f_n=f$ a.e. in $U$ whenever $n\geq\|f\|_{L^\infty(U)}$.
As in \cite[Theorem 1.2]{L25jm}, let $u_n \in H^{2,\frac{d}{2}}(U) \cap H^{1,2}_0(U) \cap L^{\infty}(U)$ denote the unique strong solution to \eqref{maineqpa} with $f$ replaced by $f_n$. Using integration by parts, we obtain for each $n \geq 1$,
\begin{equation} \label{approsoluti}
\int_U \langle \nabla u_n, \nabla \psi \rangle 
+ \langle \mathbf{H}, \nabla u_n \rangle \psi 
+ (c+\alpha) u_n \psi \, dx 
= \int_U f_n \psi \, dx
\quad \text{ for all } \psi \in C_0^\infty(U).
\end{equation}
Let $\rho$ be the function constructed in Theorem \ref{existrho}, and define $\mathbf{B}$ and $\mu := \rho\,dx$ as in \eqref{divfreebdef} and \eqref{rhodivfreeb}, respectively. Let $(\mathcal{E}, D(\mathcal{E}))$ be the Dirichlet form constructed in Proposition \ref{matebisdiri}, and let $(G_{\alpha})_{\alpha>0}$ be the strongly continuous sub-Markovian resolvent of contractions on $L^2(U, \mu)$ given in Proposition \ref{resolvpropo}. Then, by Proposition \ref{matebisdiri} and \cite[Chapter I, Theorem 2.8]{MR92}, $G_{\alpha} f_n \in D(\mathcal{E})= H^{1,2}_0(U)$ satisfies, 
\begin{align}
&\int_U \langle \rho \nabla G_{\alpha} f_n, \nabla \varphi \rangle 
+ \langle \rho \mathbf{B}, \nabla G_{\alpha} f_n \rangle \varphi 
+ \rho (c+\alpha) G_{\alpha} f_n \cdot \varphi \, dx  \nonumber \\
&  = \mathcal{E}(G_{\alpha}f_n, \varphi)  +\alpha \int_{U} G_{\alpha} f_n \cdot \varphi \,d\mu = \int_{U} f_n \varphi \,d\mu = \int_{U} \rho f_n \varphi\,dx, \quad \forall \varphi \in C_0^{\infty}(U) \text{ and } n \geq 1. \label{varidenresolv}
\end{align}
Thus, applying Theorem \ref{divfretransform} to \eqref{varidenresolv}, we get
\begin{equation} \label{apresolvari}
\int_U \langle \nabla G_{\alpha} f_n, \nabla \psi \rangle 
+ \langle \mathbf{H}, \nabla G_{\alpha} f_n \rangle \psi 
+ (c+\alpha) G_{\alpha} f_n \cdot \psi \, dx 
= \int_U f_n \psi \, dx,
\quad \forall\, \psi \in C_0^\infty(U) \text{ and } n \geq 1.
\end{equation}
Applying the uniqueness result (\cite[Theorem 1.1(ii)]{L25jm}) to \eqref{approsoluti} and \eqref{apresolvari}, we have $G_{\alpha} f_n = u_n$ in $U$. Moreover, by Proposition \ref{resolvpropo}, for each $\theta \in [1, \infty]$,
\begin{equation} \label{contraestimthe}
\|u_n\|_{L^{\theta}(U, \mu)} \leq \frac{1}{\alpha} \|f_n\|_{L^{\theta}(U, \mu)}.
\end{equation}
On the other hand, by \cite[Theorem 1.1(i)]{L25jm},
\[
\|u_n\|_{H^{1,2}_0(U)} \leq K_5 \|f_n\|_{L^{\frac{2d}{d+2}}(U)},
\]
where $K_5 > 0$ is independent of $n$ and $f$. Hence, by the weak compactness of closed and bounded convex subsets of $H^{1,2}_0(U)$ and applying the uniqueness result \cite[Theorem 1.1(ii)]{L25jm}, there exists a subsequence of $(u_n)_{n \geq 1}$, again denoted by $(u_n)_{n \geq 1}$, such that
\begin{equation} \label{weakh12conv}
\lim_{n \to \infty} u_n = u \quad \text{weakly in } H^{1,2}_0(U), 
\qquad
\lim_{n \to \infty} u_n = u \quad \text{strongly in } L^{2}(U).
\end{equation}
We now consider the following three cases. \\ \\
\noindent
\textsf{\underline{Case 1}:} $\theta = \infty$.  
From \eqref{contraestimthe}, there exist $v \in L^{\infty}(U, \mu)$ and a subsequence (still denoted $(u_n)_{n \geq 1}$) such that\begin{equation} \label{weakstarconver}
\lim_{n \to \infty} u_n = v \quad \text{weakly* in } L^{\infty}(U, \mu).
\end{equation}
By the lower semicontinuity of the norm under weak* convergence and \eqref{contraestimthe},
\begin{equation} \label{lowersemicon}
\|v\|_{L^{\infty}(U, \mu)} \leq \liminf_{n \to \infty} \|u_n\|_{L^{\infty}(U, \mu)}
\leq \liminf_{n \rightarrow \infty}\frac{1}{\alpha} \|f_n\|_{L^{\infty}(U, \mu)}=\frac{1}{\alpha} \|f\|_{L^{\infty}(U, \mu)}.
\end{equation}
Combining \eqref{weakh12conv}, \eqref{weakstarconver}, and \eqref{lowersemicon}, we deduce that $u = v \in H^{1,2}_0(U) \cap L^{\infty}(U)$ and
\begin{equation} \label{linftycontra}
\|u\|_{L^{\infty}(U, \mu)} \leq \frac{1}{\alpha} \|f\|_{L^{\infty}(U, \mu)}.
\end{equation}
Thus, \eqref{linftycontra} immediately implies \eqref{contrawithrele}.

\medskip
\noindent
\textsf{\underline{Case 2}:} $\theta \in [2, \infty)$.  
From \eqref{contraestimthe}, there exist $v \in L^{\theta}(U, \mu)$ and a subsequence (still denoted $(u_n)_{n \geq 1}$) such that
\begin{equation} \label{wealconverthe}
\lim_{n \to \infty} u_n = v \quad \text{weakly in } L^{\theta}(U, \mu).
\end{equation}
By the lower semicontinuity of the norm under weak convergence and \eqref{contraestimthe},
\begin{equation} \label{lowersemweak}
\|v\|_{L^{\theta}(U, \mu)} \leq \liminf_{n \to \infty} \|u_n\|_{L^{\theta}(U, \mu)}
\leq \liminf_{n \rightarrow \infty}\frac{1}{\alpha} \|f_n\|_{L^{\theta}(U, \mu)}=\frac{1}{\alpha} \|f\|_{L^{\theta}(U, \mu)}.
\end{equation}
Combining \eqref{weakh12conv}, \eqref{wealconverthe}, and \eqref{lowersemweak} we conclude that $u = v \in H^{1,2}_0(U) \cap L^{\theta}(U)$ and
\[
\|u\|_{L^{\theta}(U, \mu)} \leq \frac{1}{\alpha} \|f\|_{L^{\theta}(U, \mu)}.
\]
Thus, \eqref{contrawithrele} follows. \\ \\
\noindent
\textsf{\underline{Case 3}:} $\theta \in [1, 2)$.  
From \eqref{weakh12conv} and \eqref{contraestimthe},
\[
\|u\|_{L^{\theta}(U, \mu)}
= \lim_{n \to \infty} \|u_n\|_{L^{\theta}(U, \mu)} \leq \lim_{n \rightarrow \infty}\frac{1}{\alpha} \|f_n\|_{L^{\theta}(U, \mu)}= \frac{1}{\alpha} \|f\|_{L^{\theta}(U, \mu)}.
\]
\medskip
By Cases 1-3, \eqref{weightedcontra} follows.  
Finally, \eqref{weightedcontra} together with Theorem \ref{existrho}(ii) yields, for each $\theta \in [1, \infty)$,
\begin{align*}
\int_{U} |u|^{\theta}\, dx
&\leq \frac{1}{\min_{\overline{U}} \rho} \int_{U} |u|^{\theta}\, d\mu \leq \frac{1}{\min_{\overline{U}} \rho} \cdot \frac{1}{\alpha^{\theta}} \int_{U} |f|^{\theta} \, d\mu \\
&\leq \frac{\max_{\overline{U}} \rho}{\min_{\overline{U}} \rho} \cdot \frac{1}{\alpha^{\theta}} \int_{U} |f|^{\theta} \, dx \leq \frac{K_1}{\alpha^{\theta}} \int_{U} |f|^{\theta} \, dx,
\end{align*}
which implies \eqref{contrawithrele}. \\ 
\text{} \hfill $\square$
\centerline{}
{\bf Proof of Theorem \ref{pinnapplica}}\\
Fix a constant $\alpha > 0$. Let $u \in H^{2, \frac{d}{2}}(U) \cap H^{1,2}_0(U) \cap L^{\infty}(U)$ be the unique strong solution to \eqref{maineqpa} as in Theorem \ref{mainexisthm}.
Then $u-\Psi \in H^{2, \frac{d}{2} \wedge 2}(U) \cap H^{1,2}_0(U) \cap L^{\infty}(U)$ is the unique weak and strong solution to  
\begin{equation*} 
\left\{
\begin{alignedat}{2}
- \Delta (u-\Psi) + \langle \mathbf{H}, \nabla (u-\Psi) \rangle + (c+\alpha) (u-\Psi) &= f-\mathcal{L} \Psi &&\quad \text{in } U, \\ 
u-\Psi &= 0 &&\quad \text{on } \partial U,
\end{alignedat}
\right.
\end{equation*}
Observe that  
\[
f-\mathcal{L}\Psi = f + \Delta \Psi - \langle \mathbf{H}, \nabla \Psi \rangle - (c+\alpha) \Psi \in L^2(U).
\]  
Let $\rho$ be the function obtained in Theorem \ref{existrho}, and define $\mathbf{B}$ and $\mu := \rho\,dx$ as in \eqref{divfreebdef} and \eqref{rhodivfreeb}, respectively.  
Let $(\mathcal{E}, D(\mathcal{E}))$ denote the Dirichlet form given in Proposition \ref{matebisdiri}, and let $(G_{\alpha})_{\alpha>0}$ be the strongly continuous sub-Markovian resolvent of contractions on $L^2(U, \mu)$ described in Proposition \ref{resolvpropo}.  
By the $L^2(U, \mu)$ contraction property of $(G_{\alpha})_{\alpha>0}$, we have  
\begin{equation} \label{elltwocontracprop}
\| G_{\alpha} g \|_{L^2(U, \mu)} \leq \frac{1}{\alpha} \|g\|_{L^2(U, \mu)} \quad \text{ for all } g \in L^2(U, \mu).
\end{equation}
Meanwhile, by Proposition \ref{matebisdiri} and \cite[Chapter I, Theorem 2.8]{MR92},  
$G_{\alpha} (f - \mathcal{L} \Psi) \in D(\mathcal{E}) = H^{1,2}_0(U)$ satisfies  
\begin{align}
&\int_U \langle \rho \nabla G_{\alpha}(f- \mathcal{L} \Psi) , \nabla \varphi \rangle 
+ \langle \rho \mathbf{B}, \nabla G_{\alpha}(f- \mathcal{L} \Psi) \rangle \varphi 
+ \rho (c+\alpha) G_{\alpha}(f- \mathcal{L} \Psi) \cdot \varphi \, dx  \nonumber \\
& = \mathcal{E}\big(G_{\alpha}(f- \mathcal{L} \Psi), \varphi\big) +\alpha \int_{U} G_{\alpha}(f- \mathcal{L} \Psi) \cdot \varphi \,d\mu \nonumber \\
&= \int_{U} (f-\mathcal{L} \Psi) \varphi \,d\mu = \int_{U} \rho (f-\mathcal{L} \Psi)\varphi\,dx \quad \text{ for all } \varphi \in C_0^{\infty}(U). \label{varidenrephys}
\end{align}
Thus, applying Theorem \ref{divfretransform} to \eqref{varidenrephys}, we get
\begin{align*}
&\int_U \langle \nabla G_{\alpha} (f-\mathcal{L} \Psi), \nabla \psi \rangle 
+ \langle \mathbf{H}, \nabla G_{\alpha} (f-\mathcal{L} \Psi) \rangle \psi 
+ (c+\alpha) G_{\alpha} (f-\mathcal{L} \Psi) \cdot \psi \, dx  \\
&\qquad = \int_U (f-\mathcal{L} \Psi)\, \psi \, dx \quad \text{ for all } \psi \in C_0^\infty(U).
\end{align*}
Applying the uniqueness result \cite[Theorem 1.1(ii)]{L25jm}, we obtain that  
\[
G_{\alpha} (f- \mathcal{L} \Psi) = u-\Psi \quad \text{ in } U.
\]  
Finally, by \eqref{elltwocontracprop}, we get \eqref{pinnesti}.   In conclusion, \eqref{pinnesti} together with Theorem \ref{existrho}(ii) yields  
\begin{align*}
\int_{U} |u-\Psi|^{2}\, dx
&\leq \frac{1}{\min_{\overline{U}} \rho} \int_{U} |u-\Psi|^{2}\, d\mu \\
&\leq \frac{1}{\min_{\overline{U}} \rho} \cdot \frac{1}{\alpha^{2}} \int_{U} |f- \mathcal{L}\Psi|^{2} \, d\mu \\
&\leq \frac{\max_{\overline{U}} \rho}{\min_{\overline{U}} \rho} \cdot \frac{1}{\alpha^{2}} \int_{U} |f-\mathcal{L}\Psi|^{2} \, dx \\
&\leq \frac{K_1}{\alpha^{2}} \int_{U} |f-\mathcal{L} \Psi|^{2} \, dx,
\end{align*}
which proves \eqref{pracpinestim}.\\
\text{} \hfill $\square$

\begin{rem} \label{finalrema}
We note that the constant $K_1^{\frac{1}{\theta}}$ in the estimate \eqref{contrawithrele} of Theorem \ref{mainexisthm} is improved compared to the constant $K_1$ in \cite[(7) in Theorem 1.1(iii)]{L25jm}. Theorem \ref{pinnapplica} can be effectively applied to error estimates for physics-informed neural networks (PINNs).  
Indeed, in the one-dimensional case, \cite{YL24} systematically develops the $L^2$-contraction estimates and applies them to error bounds for PINNs. Suppose, a priori, that the unique solution $u$ to \eqref{maineqpa} given by \cite[Theorem 1.2]{L25jm} is continuous on $\overline{U}$.  
Under the assumptions of Theorem \ref{pinnapplica}, let $\Psi$ be
a smooth function constructed via a neural network, with zero boundary values on $\partial U$, and fixed independently of the Monte Carlo samples.
Let $(X_i)_{i \geq 1}$ be a sequence of independent and identically distributed random variables on a probability space $(\Omega, \mathcal{F}, \mathbb{P})$, uniformly distributed on $U$.  
Then, by Monte Carlo integration (see \cite[Proposition 2.4]{YL24}), the norms
\[
\| u - \Psi \|_{L^2(U)}^2 
\quad \text{and} \quad
 \, \| f - \mathcal{L} \Psi \|_{L^2(U)}^2
\]
can be approximated by
\[
\frac{|U|}{n} \sum_{i=1}^n |u(X_i) - \Psi(X_i)|^2
\quad \text{and} \quad
\frac{|U|}{ n} \sum_{i=1}^n |f(X_i) - \mathcal{L} \Psi(X_i)|^2,
\]
respectively, where $\mathcal{L}\Psi$ is defined as in \eqref{lpsidefined}.
Hence, the estimate \eqref{pracpinestim} can be reformulated in the Monte Carlo setting as
\[
 \sum_{i=1}^n |u(X_i) - \Psi(X_i)|^2
\ \leq \frac{K_1}{\alpha^2}\cdot 
\sum_{i=1}^n |f(X_i) - \mathcal{L} \Psi(X_i)|^2
\quad \text{with high probability as } n\to\infty.
\]
Therefore, as $\alpha > 0$ increases, the error-to-loss ratio, given by $\frac{K_1}{\alpha^2}$, decays to zero, implying that the error estimates remain robust for large $\alpha$ (cf. \cite{YL24}).  However, for small values of $\alpha$, the estimates may no longer be efficient, and further research is needed to analyze this setting.
\end{rem}

\section{Conclusions and discussion}
In this paper, we established improved contraction estimates for strong solutions to a class of linear elliptic diffusion-advection-reaction equations with singular coefficients. The key to our approach was the construction of a Dirichlet form and the representation of the solution as the image of the associated weighted resolvent. This method provided sharper bounds compared to previous results in \cite{L25jm}. \\
Our main contribution in Theorem \ref{mainexisthm} is the refinement of the $L^\theta$-contraction estimate, where the constant was improved from $\frac{K_1}{\alpha}$ to $\frac{K_1^{\frac{1}{\theta}}}{\alpha}$. This improvement is particularly significant for large values of $\theta$. Notably, for the crucial case of the $L^\infty$-contraction, we derived an explicit constant of $\frac{1}{\alpha}$, which depends only on the parameter $\alpha$ and is independent of the domain geometry or the norms of the lower-order coefficients. This result offers a clearer and more precise quantitative understanding of the solution's behavior.\\
In the context of scientific machine learning, these results provide a rigorous error analysis for Physics-Informed Neural Networks (PINNs). As demonstrated in Theorem \ref{pinnapplica} and discussed in Remark \ref{finalrema}, the derived estimates provide a rigorous error bound for PINN approximations. The error estimate in Theorem~\ref{pinnapplica}
is robust with respect to large $\alpha$,
since $K_1$ is independent of $\alpha$ and
the constant $\frac{\sqrt{K_1}}{\alpha}$ decreases
as $\alpha$ increases.\\
However, our analysis also highlights a limitation and a clear direction for future work. The derived error bounds get worse as $\alpha$ approaches zero, indicating that these contraction estimates may lose their effectiveness for equations with weak reaction terms. Therefore, subsequent research should focus on developing new analytical tools that work well even when $\alpha$ is small. \\
A promising direction for future research involves extending our framework to accommodate more general diffusion coefficients $A=(a_{ij})_{1 \leq i,j \leq d}$. A key objective would be to remove the vanishing mean oscillation ($VMO$) assumption, relying only on the standard uniform ellipticity condition. Furthermore, weakening the conditions on the zero-order term, for instance, by allowing $c \in L^1(U)$ instead of $c \in L^{\frac{d}{2}}(U)$ and relaxing its positivity, presents another valuable avenue of research. \\ \\

\section*{Acknowledgments}
This work was supported by the National Research Foundation of Korea (NRF) grant funded by the Korea government (MSIT) (RS-2026-25598810) and the Gyeongbuk ANCHOR system-(Regional Growth Innovation LAB) program through the Gyeongbuk ANCHOR Center, funded by the Ministry of Education (MOE) and the Gyeongsangbuk-do, Republic of Korea (2026-anchor-15-105).

\centerline{}
\centerline{}
Haesung Lee\\
Department of Mathematics and Big Data Science,  \\
Kumoh National Institute of Technology, \\
Gumi, Gyeongsangbuk-do 39177, Republic of Korea, \\
E-mail: fthslt@kumoh.ac.kr, \; fthslt14@gmail.com
\end{document}